\documentclass{amsart}

\usepackage{amsmath,amsfonts,amsthm,amssymb,color}

\newtheorem{theorem}[equation]{Theorem}

\newtheorem{lemma}[equation]{Lemma}

\newtheorem{corollary}[equation]{Corollary}

\theoremstyle{remark}

\theoremstyle{definition}
\newtheorem{defn}[equation]{Definition}

\newtheorem*{prop22}{Proposition 2.2 (Joshi-S\`a Barreto)}
\newtheorem*{thm31}{Theorem 3.1 (Joshi-S\`a Barreto)}
\newtheorem*{prop31}{Proposition 3.1 (Joshi-S\`a Barreto)}
\newtheorem*{thm41}{Theorem 4.1 (Joshi-S\`a Barreto)}
\newtheorem*{thm42}{Theorem 4.2 (Joshi-S\`a Barreto)}

\numberwithin{equation}{section}

\renewcommand{\qed}{\hspace*{\fill} \setlength{\unitlength}{1mm}
\begin{picture}(2.5,2.5)
      \put(0,0){\framebox(2.5,2.5){}}
  \end{picture}
\setlength{\unitlength}{1pt}}

\newcommand{\reals}{{\bf R}}

\newcommand{\tr}{{\rm tr}}

\newcommand{\cD}{{\mathcal{D}}}
\newcommand{\cN}{{\mathcal{N}}}
\newcommand{\cP}{{\mathcal{P}}}
\newcommand{\cR}{{\mathcal{R}}}

\newcommand{\cV}{\mathcal{V}}
\newcommand{\cO}{\mathcal{O}}
\newcommand{\cF}{\mathcal{F}}

\newcommand{\N}{\mathbb{N}}
\newcommand{\R}{\mathbb{R}}
\newcommand{\C}{\mathbb{C}}

\newcommand{\bbH}{\mathbb{H}}

\newcommand{\cC}{\mathcal{C}}
\newcommand{\cL}{\mathcal{L}}
\newcommand{\fp}{\mathfrak{p}}

\begin{document}
\title[Dynamics of Asymptotically Hyperbolic Manifolds]{Dynamics of Asymptotically Hyperbolic Manifolds}

\author[Julie Rowlett]{Julie Rowlett}
\address{Department of Mathematics, South Hall 6607, University of California, Santa Barbara, CA 93106.} \email{rowlett@math.ucsb.edu}


\begin{abstract}

We prove a dynamical wave trace formula for asymptotically hyperbolic $(n+1)$ dimensional manifolds with negative (but not necessarily constant) sectional curvatures which equates the renormalized wave trace to the lengths of closed geodesics.  This result generalizes the classical theorem of Duistermaat-Guillemin \cite{DG} for compact manifolds and the results of \cite{GN}, \cite{GZ}, and \cite{P} for hyperbolic manifolds with infinite volume.  A corollary of this dynamical trace formula is a dynamical resonance-wave trace formula for compact perturbations of convex co-compact hyperbolic manifolds which we use to prove a growth estimate for the length spectrum counting function. We next define a dynamical zeta function and prove its analyticity in a half plane.  In our main result, we produce a prime orbit theorem for the geodesic flow.  This is the first such result for manifolds which have neither constant curvature nor finite volume.  As a corollary to the prime orbit theorem, using our dynamical resonance-wave trace formula, we show that the existence of pure point spectrum for the Laplacian on negatively curved compact perturbations of convex co-compact hyperbolic manifolds is related to the dynamics of the geodesic flow.   

\end{abstract}
\maketitle


\section{Introduction}

The spectral theory and dynamics of hyperbolic manifolds has interested mathematicians for many years.  Motivated by recent developments in theoretical physics, mathematical interest in \em asymptotically hyperbolic \em manifolds has been piqued; see for example \cite{MMAH}, \cite{BP}, \cite{GZ2}, \cite{JSB2}, \cite{P}, \cite{GrZ}, and  \cite{FG}.  Asymptotically hyperbolic manifolds are appealing because they arise in connections with conformal field theory--anti-de Sitter correspondence, \cite{Gr} and \cite{SbZ}, and because they are a class of manifolds on which geometric scattering theory can be developed; see \cite{JSB}, \cite{G0}, and \cite{G1}.  

Recall that on a compact manifold, the spectrum of the Laplacian $\Delta$ is a discrete subset of $\R^+,$ 
$$\sigma(\Delta) = \{ \lambda_k ^2 \}_{k=1} ^{\infty},$$
and the wave trace is formally the distribution 
$$\sum_{k \geq 1} e^{i \lambda_k t}.$$
The singularities of the wave trace occur at $t=0$ and at the lengths of closed geodesics.  The principle part of the singularities for $t>0$ was computed in the seminal work of Duistermaat and Guillemin \cite{DG} to be 
$$ \frac{l(\gamma)}{\sqrt{|\det( I - P_{\gamma} ^k)|}} \delta(|t| - k l(\gamma)),$$
who also computed the asymptotics of the ``big singularity'' at $t=0.$  One then has at least formally,
$$\sum_{k \geq 1} e^{i \lambda_k |t|} = \sum_{ \{\gamma \} } \sum_{k=1} ^{\infty} \frac{l(\gamma)  \delta(|t| - k l(\gamma))}{\sqrt{|\det( I - P_{\gamma} ^k )|}} + A(t),$$
where $\{ \gamma \}$ are the primitive closed geodesics with length $l(\gamma),$ $P^k _{\gamma}$ is the $k$-times Poincar\'e map in the cotangent bundle about $\gamma,$ and the remainder $A(t)$ is exponentially singular at $0$ and smooth for $t>0.$  This formula shows a beautiful connection via the wave trace between the Laplace and length spectra.  

On compact hyperbolic surfaces this dynamical formula is well known.  The remainder term $A(t)$ is explicitly computable as a ratio of hyperbolic trigonometric functions and in this case, the result is known as a Selberg trace formula and also holds for hyperbolic surfaces with infinite volume and cusps \cite{GZ}.  For cofinite surfaces, congruence subgroups, and $PSL(2, \R),$ Hejhal has proven Selberg trace formulae, \cite{Hejc}, \cite{Hej}, \cite{Hej2}.  In higher dimensions, Gangolli and Warner have proven Selberg trace formulae, \cite{GW}, \cite{GW2}.  Further results include \cite{art}, \cite{bjp}, \cite{pp}, \cite{juhl}, and \cite{mu}.  On non-compact manifolds with constant negative curvature, algebraic group averaging methods have been used to prove the trace formula \cite{GZ}, \cite{GN}, and \cite{P}.  Since these hyperbolic manifolds may have infinite volume, a renormalized wave trace replaces the standard wave trace in the dynamical formula.  This renormalized wave trace is known as the $0$-trace and was introduced by Guillop\'e and Zworski in \cite{GZ4}.  Our dynamical wave trace formula generalizes \cite{DG} to asymptotically hyperbolic manifolds.   

\begin{theorem}  \label{theorem1} 
Let $(X, g)$ be an asymptotically hyperbolic $n+1$ dimensional manifold with negative sectional curvatures. Let 0-tr $\cos \left( t \sqrt{ \Delta - \frac{n^2}{4}} \right)$ denote the regularized trace of the wave group, (\ref{0trcos}).  Let $\cL_p$ denote the set of primitive closed geodesics of $(X, g),$ and for $\gamma \in \cL_p,$ let $l(\gamma)$ denote the length of $\gamma.$  Then, for any $\epsilon >0,$
\begin{equation} \label{trace1} \textrm{0-tr} \cos \left( t \sqrt{ \Delta - \frac{n^2}{4}} \right) = \sum_{\gamma \in \cL_p, \, k \in \N} \frac{ l(\gamma) \delta(|t| - k l(\gamma))}{ \sqrt{ | \det(I - \cP^k _{\gamma}) | }} + A(t),\end{equation}
as a distributional equality in $\cD'([\epsilon, \infty)),$\footnote{$\cD'(X)$ is the dual of $\cC^{\infty} _0 (X).$} where the remainder $A(t)$ is smooth and bounded as $t \to \infty,$ and $\cP^k _{\gamma}$ is the $k$-times Poincar\'e map around $\gamma$ in the cotangent bundle.  
 \end{theorem}

As a corollary to this theorem we prove a dynamical resonance-wave trace formula.  

\begin{corollary} \label{corollary1} 
Let $(X, g)$ be a negatively curved\footnote{We will use ``negatively curved'' to mean all sectional curvatures are negative.} compact perturbation of a convex cocompact $n+1$ dimension hyperbolic manifold, (\ref{pert}).  Then, for any $\epsilon>0,$ as an element of $\cD'([\epsilon, \infty)),$ we have the distributional equality  
$$\textrm{0-tr} \cos \left( t \sqrt{ \Delta - \frac{n^2}{4}} \right) = \frac{1}{2} \sum_{s \in \cR^{sc} } e^{(s - n/2)|t|} + B(t) = \sum_{\gamma \in \cL_p, \, k \in \N} \frac{ l(\gamma) \delta(|t| - k l(\gamma))}{ \sqrt{ | \det(I - \cP^k _{\gamma})| }} + C(t).$$
Above, $B(t)$ and $C(t)$ satisfy the regularity and bounds for $A(t)$ in (\ref{trace1}), and we note that the resonances $\cR^{sc}$ of the scattering operator (\ref{scat})  are summed \em with multiplicity,  \em where the actual spectral parameter is $\Lambda = s(n-s).$ 
\end{corollary}

Recall that the Laplacian on an asymptotically hyperbolic manifold of dimension $n+1$ has absolutely continuous spectrum, $\sigma_{ac} (\Delta) = [ \frac{n^2}{4}, \infty ),$ and a finite pure point spectrum, $\sigma_{pp} (\Delta) \subset (0, \frac{n^2}{4} ).$  As a corollary to the dynamical resonance-wave trace formula, we prove the following preliminary growth estimate for the length spectrum counting function.

\begin{corollary} \label{corollary2}
Let $(X, g)$ be a negatively curved compact perturbation of a convex cocompact $n+1$ dimension hyperbolic manifold.  
Assume that $\sigma_{pp}(\Delta) \neq \emptyset,$ and let 
$$\Lambda_0 = \inf \{ \Lambda \in \sigma_{pp} (\Delta) \},$$
with corresponding $s_0$ satisfying $\Lambda_0 = s_0 (n-s_0).$  Let $\cL$ be the set of closed geodesics in $X,$ and for $\gamma \in \cL,$ let $l(\gamma)$ be the length of $\gamma.$  Then, the length spectrum counting function
\begin{equation} \label{eq:lect} N(T) := \# \{ \gamma \in \cL : l(\gamma) \leq T \} \end{equation} 
satisfies
$$\limsup_{T \to \infty} \frac{T N(T)}{e^{(s_0 + (k_1 - 1)n/2)T}} \leq 1 \leq \liminf_{T \to \infty} \frac{T N(T)}{e^{(s_0 + (k_2 - 1)n/2)T}},$$
where $0 < k_2 \leq 1 \leq k_1,$ and the sectional curvatures $\kappa$ satisfy $-k_1 ^2 \leq \kappa \leq -k_2 ^2$ on $X.$  
 \end{corollary}
Although this result is not sharp like the prime orbit theorem, it is appealing for its proof which demonstrates a key difference between Euclidean and hyperbolic scattering:  in hyperbolic scattering it is possible to have exponential growth in the spectral side of the trace formula, whereas in Euclidean scattering the spectral side has purely oscillatory terms.  The next result is analyticity of the dynamical zeta function in a half plane determined by the topological entropy of the geodesic flow.

\begin{theorem}  \label{theorem2}
Let $(X, g)$ be an asymptotically hyperbolic negatively curved $n+1$ dimensional manifold.  Let $L_p$ be the set of primitive closed orbits of the geodesic flow on $X,$ and for $\gamma \in L_p,$ let $l_p (\gamma)$ be the length of the primitive period of $\gamma.$ Then, the dynamical zeta function 
$$ Z(s) = \exp \left( \sum_{\gamma \in L_p} \sum_{k \in \N} \frac{e^{-k s l_p(\gamma)}}{k} \right)$$
converges absolutely for $\mathfrak{R}(s) > h,$ where $h$ is the topological entropy of the geodesic flow (\ref{topent}).  The weighted dynamical zeta function,
$$\tilde{Z} (s) = \exp \left( \sum_{\gamma \in L_p} \sum_{k \in \N} \frac{e^{-k s l_p(\gamma)}}{k\sqrt{ | \det(I - \cP^k _{\gamma})|}} \right),$$  
converges absolutely for $\mathfrak{R}(s) > \fp(-H/2),$ where $\fp$ is the topological pressure, (\ref{press}), and $H$ is the Sinai-Ruelle-Bowen potential, (\ref{srb}). \end{theorem}

Our main result produces a prime orbit theorem for the geodesic flow on negatively curved $n+1$ dimension asymptotically hyperbolic manifolds.  

\begin{theorem} \label{theorem3} Let $(X, g)$ be an asymptotically hyperbolic $n+1$ dimensional manifold with negative sectional curvatures.  Let $L_p$ be the set of primitive closed orbits of the geodesic flow, and for $\gamma \in L_p,$ let $l_p(\gamma)$ be the length of the primitive period of $\gamma.$  Let $h$ be the topological entropy of the geodesic flow, (\ref{topent}).  The dynamical zeta function, 
$$Z(s) = \exp \left( \sum_{\gamma \in L_p} \sum_{k \in \N} \frac{e^{-k s l_p(\gamma)}}{k} \right)$$
has a nowhere vanishing analytic extension to an open neighborhood of $\mathfrak{R}(s) \geq h$ except for a simple pole at $s=h.$  Moreover, the length spectrum counting function (\ref{eq:lect}) satisfies 
$$\lim_{T \to \infty} \frac{T N(T)}{e^{hT}} = 1.$$ 
\end{theorem} 

Finally, we use the prime orbit theorem and the trace formula Corollary \ref{corollary2} to prove a result which shows that the existence of pure point spectrum depends on the topological entropy of the geodesic flow and the curvature bounds for negatively curved compact perturbations of convex co-compact hyperbolic manifolds.  

\begin{corollary} \label{corollary3}
Let $(X, g)$ be a negatively curved compact perturbation of a convex cocompact $n+1$ dimension hyperbolic manifold with topological entropy $h.$  Then there exist $0 < k_2 \leq 1 \leq k_1$ so that the sectional curvatures $\kappa$ satisfy $-k_1 ^2 \leq \kappa \leq -k_2 ^2.$  If the topological entropy $h > \frac{nk_1}{2},$ then $\sigma_{pp} (\Delta) \neq \emptyset,$ and moreover, there is $\Lambda_0 = s_0 (n-s_0) \in \sigma_{pp} (\Delta)$ with $s_0 \geq h + \frac{(n-k_1)}{2}.$  Conversely, if the topological entropy $h \leq \frac{nk_2}{2},$ then $\sigma_{pp}(\Delta) = \emptyset.$  
\end{corollary}

This paper is organized as follows:  in section 2, we recall some basic spectral and geometric properties of asymptotically hyperbolic manifolds including the $0$-renormalization and the results of \cite{JSB2} for the wave group.  In section 3, we recall some basic definitions and results from dynamics and prove the dynamical trace formula and its corollaries.  In section 4, we define the dynamical zeta function and study its regularity.  We show in section 5 that we may apply Parry and Pollicott's work for Axiom A flows to our setting, and according to their work we produce a prime orbit theorem for the geodesic flow and its corollary relating the pure point spectrum to the dynamics of the flow.  Concluding remarks comprise section 6.  

The author would like to thank I. Polterovich and D. Jakobson for motivating this work and for several useful discussions and correspondence and to thank the Centre de R\'echerches Math\'ematiques in Montr\'eal, Canada where this work was initiated.  

\section{Asymptotically Hyperbolic Manifolds}

A manifold with boundary $(X^{n+1}, \partial X)$ is \em asymptotically hyperbolic \em if there exists a boundary defining function $x$ so that a neighborhood of $\partial X$ admits a product decomposition $(0, \epsilon)_x \times \partial X,$ and with respect to this decomposition the metric takes the form 
\begin{equation}\label{ahmetric} g = \frac{dx^2 + h(x, y, dx, dy)}{x^2} \end{equation}
where $h|_{\{x=0\}}$ is independent of $dx.$  We note that there is no one canonical metric on $\partial X$ but rather a conformal class of metrics induced by $h|_{\{x=0\}} = h_0,$ and $(\partial X, [h_0])$ is called the conformal infinity; see for example \cite{GrZ} and \cite{G1}.  

It was observed by Mazzeo and Melrose in \cite{MMAH} that $X$ is a complete Riemannian manifold; moreover, along any smooth curve in $X - \partial X$ approaching a point $p \in \partial X,$ the sectional curvatures of $g$ approach $- |d x|^2 _{x^2 g}.$  For each $h \in [h_0]$ there exists a unique (near the boundary) boundary defining function $x$ so that 
$$|dx|_{x^2 g} = 1,$$
near $Y := \partial X.$  With this normalization, the sectional curvatures approach $-1$ at $\partial X,$ hence the name ``asymptotically hyperbolic.'' 

A large class of interesting asymptotically hyperbolic metrics are the conformally compact metrics.  In particular, four dimensional conformally compact Einstein metrics have recently received much attention in both geometric analysis and mathematical physics due to their relation to quantum field theory and quantum gravity \cite{albin}, \cite{FG}.  Of course, the simplest example of an asymptotically hyperbolic metric is the hyperbolic metric 
$$\frac{dx^2 + dy^2}{x^2},$$
on the upper half space $\bbH^{n+1} = \R^+ _x \times \R^n _y.$  It is noted in \cite{JSB} that for a metric which has the form 
$$g = \frac{dx^2 + h(x, y, dx, dy)}{x^2}$$
there exists a product decomposition so that
$$g = \frac{dx^2 + h(x, y, dy)}{x^2}.$$

\subsection{The wave group}

Let $(X, g)$ be an asymptotically hyperbolic manifold of dimension $n+1.$  The Laplacian on an asymptotically hyperbolic manifold has absolutely continuous spectrum, $\sigma_{ac} (\Delta) = [ \frac{n^2}{4}, \infty ),$ and a finite pure point spectrum, $\sigma_{pp} (\Delta) \subset (0, \frac{n^2}{4} ).$  It is natural in this setting to use the spectral parameter $s(n-s).$  In \cite{MMAH}, completed by \cite{G0}, it was shown that if $g$ is even modulo $O(x^{2k+1}),$ then the resolvent $R(s) = (\Delta - s(n-s))^{-1}$ extends meromorphically from $\{ \mathfrak{R} (s) > \frac{n}{2} \}$ to $\C - 
\{ \frac{n+1}{2} - k - \N \}.$  Recall that for a product decomposition near $\partial X$ so that $g = \frac{dx^2 + h(x, y, dy)}{x^2},$ $g$ is said to be \em even modulo $O(x^{2k+1})$ \em if $h$ admits a formal power series expansion, $h \sim h_0 + x h_1 + x^2 h_2 + \ldots$ so that $h_j = 0$ for $j \leq 2k.$


The wave kernel is the Schwartz kernel of the fundamental solution to 
$$\left( \partial_t ^2 + \Delta - \frac{n^2}{4} \right)U(t, w, w') = 0, \quad U(0, w, w') = \delta(w-w'), \quad \frac{\partial}{\partial t}  U(0, w, w') = 0.$$
Due to the semi-group property with respect to time, the wave kernel is also referred to as the wave group and written
$$\cos \left( t \sqrt{ \Delta - \frac{n^2}{4} }\right).$$
In \cite{JSB2}, Joshi and S\`a Barreto constructed the wave group as an element of an operator calculus on a manifold with corners obtained by blowing up $\R^+ \times X \times X.$  We briefly recall this construction and the main results of \cite{JSB2}.  This construction was heavily influenced by Melrose's work with $b$-manifolds in \cite{tapsit}.  


\subsubsection{The 0-blowup}
This construction is from \cite{JSB2}.  Following the notation introduced in \cite{tapsit}, the \em zero blowup \em of $X \times X$ is
$$X \times_0 X := [X \times X; \Delta(Y\times Y)],$$
where $\Delta(Y \times Y)$ is the diagonal in $Y \times Y,$ and $Y = \partial X.$  Recall that for a p-submanifold\footnote{A submanifold is called a ``p-submanifold'' if it admits a consistent local product decomposition.}  which in this example is $\Delta(Y \times Y) \subset X \times X,$ the radial blowup of $X \times  X$ around $\Delta(Y \times Y)$ is defined to be the union,
$$[X \times X ; Y \times Y] := \left( \left( X \times X \right) - \Delta(Y \times Y) \right) \cup N^+ ( \Delta( Y \times Y)).$$
Above, $N^+ (\Delta(Y \times Y))$ is the inward pointing spherical normal bundle of $\Delta(Y \times Y).$  For details of this construction see \cite{tapsit} or \cite{acc}.  The notation $[ A; B]$ indicates that the p-submanifold $B \subset A$ is ``blown up'' by replacing $B$ with its inward pointing spherical normal bundle.  As a consequence, smooth functions on $A$ and polar coordinates around $B$ lift to be smooth on $[A; B].$  There is a natural ``blow down'' map $\beta_*: [A; B] \to A$ and ``blow up'' map $\beta^* : A \to [A; B],$ with respect to which 
$$[A; B] = (A - B) \cup \beta^* (B).$$ 
$X \times_0 X$ has three boundary faces; the face created by blowing up is called the ``front face'' and written $\mathcal{F},$ while the other two faces are called the ``side faces.''\footnote{In \cite{JSB2}, the side faces are called the top and bottom faces.}  Figure 1 of \cite{JSB2} is an illustration of $X \times_0 X.$  

The zero blowup is analogous to the $b$-blowup of \cite{tapsit} because it creates a natural setting in which to study geometric differential operators on a Riemannian manifold with boundary whose metric has a particular structure at the boundary.  The Laplacian on an asymptotically hyperbolic manifold is locally the product of vector fields that vanish at $\partial X,$ hence the use of ``0'' in notation.\footnote{One could refer to asymptotically hyperbolic manifolds $0$-manifolds, but this is less effective in conveying the geometry of the manifold near the boundary.}   On a $\cC^{\infty}$ manifold with boundary $X,$ the space $\cV_0 (X)$ of smooth vector fields that vanish on the boundary is a Lie algebra.  In local coordinates $(x, y_1, \ldots, y_n)$ near $\partial X,$ $\cV_0(X)$ is the space of all smooth sections of a vector bundle over $X,$
$$\cV_0 (X) = \cC^{\infty} (X, ^0 TX).$$
$\cV_0 (X)$ is locally spanned by $\{ x \partial_x, x \partial_{y^1}, \ldots, x \partial_{y^n} \}.$ 

\subsubsection{The wave kernel construction}
Coupling  Melrose's techniques which extend classical microlocal analysis on compact manifolds to manifolds with boundary together with H\"ormander's study of Fourier integral operators \cite{ho}, Joshi and S\`a Barreto constructed the wave kernel on asymptotically hyperbolic manifolds \cite{JSB2}.  They assumed the Laplacian and wave operator act on sections of the half density bundle.  The metric $g$ in a neighborhood of the boundary is
$$g = \frac{dx^2 + h(x, y, dy)}{x^2},$$ 
and $g$ induces a canonical trivialization of the 1-density bundle by taking $\theta = \sqrt{\textrm{vol}(g)} |dx dy|.$  The square root of this is a natural trivialization of the half-density bundle.  The sections of the density bundle $^0 \Omega (X)$ are defined to be smooth multiples of the Riemannian density which in local coordinates is
$$\theta = h(x,y) \frac{dx}{x}\frac{dy}{x^n}, \, h \in \cC^{\infty} (X), \, h \neq 0.$$
The bundle $^0 \Omega^{\frac{1}{2}}(X)$ is the half-density bundle obtained from $^0 \Omega (X).$  The bundle $^0 \Omega^{\frac{1}{2}} (X \times X)$ is similarly defined, and the bundle $^0 \Omega^{\frac{1}{2}}(X \times_0 X)$ is defined to be the lift of $^0 \Omega^{\frac{1}{2}} (X \times X)$ under the blow up map $\beta.$  

In \cite{JSB2}, the wave kernel is constructed as an element of a zero pseudodifferential operator calculus.  The calculus is written $\Psi^{m, a, b} _0 \left(X, ^0 \Omega^{\frac{1}{2}} (X) \right),$ where $m \in \R,$ and $a, b \in \C.$  The kernel $K$ of a zero pseudodifferential operator in $\Psi^{m, a, b} _0 \left(X, ^0 \Omega^{\frac{1}{2}} (X) \right)$ was defined in \cite{MMAH} to be a distribution which can be written as $K = K_1 + K_2,$ where the lift of $K_1$ is conormal of order $m$ to the lifted diagonal $D_0 \subset X\times_0 X,$ and smooth up to the front (blown up) face, and vanishes at the remaining (side) boundary faces.  The second term, $K_2,$ is of the form $K_2 = \rho^a (\rho')^b F,$ $F \in \cC^{\infty} \left( X \times_0 X, ^0 \Omega^{\frac{1}{2}} (X \times_0 X) \right),$ where $\rho$ and $\rho'$ are defining functions for the side faces of $X \times_0 X.$  A key observation of \cite{JSB2} is that wave kernel information propagates at finite speed, and in this setting the distance to the boundary of $X$ is infinite.  Consequently, there is no support on the side faces and only on the interior of the front face.  

For a point $p \in \partial X,$ the inward pointing tangent vectors at $p,$ $X_p \subset T_p (X),$ is a manifold with boundary and metric given by 
$$g_p = (dx)^{-2} G(p),$$
where $g = x^{-2} G.$  We see that $(X_p, g_p)$ is isometric to the hyperbolic upper half space by regarding $G(p)$ and $dx$ as linear functions on $X_p.$  The fiber over a point $p \in \Delta(\partial X \times \partial X)$ was observed by \cite{MMAH} to be naturally isomorphic to $X_p.$  Since the kernel of an element $A \in \Psi^{m, a, b} _0 (X)$ is conormal to the lift of the diagonal, it can be restricted to the fiber $F_p \cong X_p,$ and this defines the kernel of the \em normal operator \em  of $A.$  Namely, the kernel of the normal operator, $N_p (A)$ is defined by 
\begin{equation} \label{normop} k(N_p (A)) := k(A) \big|_{F_p}.  \end{equation}

To define a class of 0-Fourier integral operators, \cite{JSB2} studied the symplectic structure of $^0 T^* X$ and defined a class of $0$-canonical transformations.  Such a transformation between $\cC^{\infty}$ manifolds with boundary $X$ and $Y$ is a smooth homogeneous map $\chi$ defined on an open conic set $U \subset \, \,  ^0 T^*X$ and maps $^0 T^*_{\partial X} X $ to $^0 T^* _{\partial Y} Y.$  A homogeneous canonical transformation $\chi: \, \, ^0 T^* X \to \, \, ^0 T^* X$ whose projection onto the base space is the identity when restricted to $\partial X$ is called a \em liftable canonical transformation.  \em  Certain $0$-canonical transformations from $X$ to $X$ define Lagrangian submanifolds of $T^* (X \times_0 X).$  In particular, \cite{JSB2} proved the following.

\begin{prop22} \em Let $p \in \cC^{\infty} (^0 T^* X)$ be the length function for the metric $g$ on $X,$ and let $^0 H_p$ be the 0-Hamiltonian vector field of $p.$  Let $\tilde{p}$ be the lift of $p$ to $T^* (X \times_0 X).$  For $s >0,$ let $\chi_s : \, \, ^0 T^* X \to \, \, ^0 T^* X$ be the map\mbox{~} defined by 
$$\chi_s (q) := \exp (s \, ^0 H_p )(q).$$
Then the graph of $\chi_s$ defines a smooth extendible Lagrangian submanifold of  $T^* (X \times_0 X).$  Moreover, the intersection 
$$\Lambda_s \cap T^* _{\mathcal{F}} (X \times_0 X)$$
is a smooth Lagrangian submanifold of  $T^* \mathcal{F}$ which is given by $\exp ( s H_{p_0} )(T^* _{\Delta(X \times X) \cap \mathcal{F}} \mathcal{F}),$ where $p_0 = \tilde{p} \big|_{\mathcal{F}}.$ \em 
\end{prop22}
Recall that a smooth conic closed Lagrangian submanifold $\Lambda \subset T^*(X \times _0 X)$ is \em extendible \em if it intersects $T^* _{\mathcal{F}}$ transversally, and in that case there exists a smooth conic closed Lagrangian submanifold $\Lambda_e \subset T^* (X \times_0 X)_d$ such that
$$\Lambda = \Lambda_e \cap T^* (X \times_0 X), \, \Lambda_0 = \Lambda \pitchfork T^*_{\cF} ( X \times_0 X).$$
Since the kernel of the wave group is a distribution in $\R \times X \times X,$ \cite{JSB2} defined a larger class of Lagrangians using the bundle $T^* \R \times ^0 T^* X \times ^0 T^* X$ and the canonical 1-form on $^0 T^* X \times ^0 T^*X,$ 
$$\alpha = \tau dt + \frac{\lambda}{x} dx + \frac{\mu}{x} \cdot dy - \frac{\lambda'}{x'} dx' + \frac{\mu'}{x'} \cdot dy'.$$  Let 
\begin{equation} \label{c} C = \{ (t, \tau, x, y, \lambda, \mu, x', y', \lambda', \mu') : \tau + \sqrt{p(x,y,\lambda, \mu)} = 0; (x', y', \lambda', \mu') = \chi_t (x, y, \lambda, \mu) \}, \end{equation}
where $p$ and $\chi_t$ are defined in Proposition 2.2.  The space of 0-Fourier integral operators associated to a liftable 0-canonical transformation $\chi$ is defined by \cite{JSB2} to consist of the kernels
$$I^{m,s} _0 (X, \chi, ^0 \Omega^{\frac{1}{2}} ) = \{ K \in I^{m,s} \left( X \times_0 X, \Lambda_{\chi}, ^0 \Omega^{\frac{1}{2}} \right);$$
$$K \textrm{ \em vanishes in a neighborhood of \em} \partial( X \times_0 X) - \cF \}.$$
Recall that for an extendible Lagrangian $\Lambda \subset T^* (X \times_0 X)$ or $\Lambda \subset T^* \R \times T^* (X \times_0 X),$ $I^{m,s} (\Lambda)$ is defined to be $R^s I^m (\Lambda),$ where $R$ is a defining function for the front face.  The symbol of $A$ at the front face is defined to be the restriction of $R^{-s} A$ to $\cF.$  This depends on $R$ but is invariant as a section of the normal bundle raised to the $s.$  For $C$ as in (\ref{c}), \cite{JSB2} showed that a Lagrangian submanifold $\Lambda_C \subset T^* \times T^* (X \times_0 X)$ is given by 
$$\Lambda_C := \{ (t, \tau, Y, G) : \tau + \sqrt{\tilde{p}(Y, G)} = 0, \, (Y, G) \in \Lambda_t \}.$$
Then, \cite{JSB2} defined
\begin{center}
$$I^{m,s} _0 (\R \times X, X; C, ^0 \Omega^{\frac{1}{2}}) = \{ K \in I^{m,s} \left(R \times X \times_0 X, \Lambda_C, ^0 \Omega^{\frac{1}{2}} \right);$$
$$ K \textrm{ \em vanishes in a neighborhood of \em} \partial (\R \times X \times_0 X) - (\R \times \cF) \}.$$
\end{center}

\begin{prop31} \em 
The normal operator (\ref{normop}) defines an exact sequence
$$0 \longrightarrow I^{m,1} _0 (\R \times X, X; C, ^0 \Omega^{\frac{1}{2}} ) \longrightarrow I^{m,0} _0 (\R \times X, X; C, ^0 \Omega^{\frac{1}{2}} ) \longrightarrow I^m (\cF, \Lambda^0 _C, \Omega^{\frac{1}{2}} ).$$ 
\em \end{prop31}

As a consequence of the exact sequence defined by the normal operators, the main result of \cite{JSB2} on the wave group is the following.  

\begin{thm31} \em For $t \in \R,$ let $C$ be defined by (\ref{c}).  The wave group $U(t)$ satisfies
$$U(t) = \cos\left( t \sqrt{ \Delta - \frac{n^2}{4}} \right) \in I^{-\frac{1}{4}, 0} _0 \left(\R \times X, X; C, ^0 \Omega^{\frac{1}{2}} \right).$$ \em 
\end{thm31}

The proof is very similar to the construction in section 7 of \cite{tapsit}.  After computing the normal operator, it follows from the usual theory of Fourier integral operators that there is an approximate kernel for $\cF,$ $U_0 (t) \in I^{-\frac{1}{4}} \left( \R \times \cF , \Lambda^0 _C, \Omega^{\frac{1}{2}} \right).$  By surjectivity of the exact sequence defined by the normal operator, one can pick an element $u_0 \in I^{-\frac{1}{4}, 0} _0 ( \R \times X, X; C, ^0 \Omega^{\frac{1}{2}} ).$  A similar construction applies to the second copy of $X,$ and this process is iterated to give a kernel with error in $I^{-\frac{1}{4}, \infty} _0 (\R \times X, X; C, ^0 \Omega^{\frac{1}{2}}),$ and an error in the Cauchy data which vanishes to infinite order at $\cF,$ and which is a pseudodifferential operator of order zero.  The error term is extended to be identically zero across the front face and removed in the standard way using H\"ormander's Lagrangian calculus.  We note that this iterative construction is a necessary and key ingredient for calculating the renormalized wave trace.  

\subsection{0-Renormalized integrals and the renormalized wave trace}

Like hyperbolic space, asymptotically hyperbolic spaces have infinite volume; to take traces one must introduce an integral renormalization.  Recall that the finite part $FP_{\epsilon \to 0} f(\epsilon)$ is defined as $f_0$ when $f(\epsilon) = f_0 + \sum_k f_k \epsilon^{- \lambda_k} (\log{\epsilon})^{m_k} + o(1),$ with $\mathfrak{R}( \lambda_k ) \geq 0,$ $m_k \in \N \cup \{0\}.$  Then $f_0$ is unique as shown for example in \cite{Ho1}.  The following definition is due to Guillop\'e and Zworski \cite{GZ2}.  
\begin{defn} \em
The $0$-regularized integral $\int \limits^{0 \mbox{~}\mbox{~}\mbox{~}} \omega$ of a smooth function (or density) $f$ on $X$ is defined, if it exists, as the finite part
$$\int \limits^{0 \mbox{~}\mbox{~}\mbox{~}} _X  f := FP_{\epsilon \to 0} \int_{x(p) > \epsilon} f(p) \textrm{dvol}_g(p),$$
where $x$ is a boundary defining function.  \em
\end{defn}
This naturally then depends on the boundary defining function $x,$ except in certain cases.  For example, it was shown in \cite{albin} and \cite{Gr} that the 0-volume of $X$ defined to be the 0-integral of the constant function $1$ is independent of the choice of $x$ if the dimension $n+1$ of $X$ is even.  For an operator $A$ with smooth Schwartz kernel $A(x, y)$ on $X \times X,$ we may then define the 0-trace of $A$ to be
$$0-tr (A) :=  \int \limits^{0 \mbox{~}\mbox{~}\mbox{~}} _X A(x, x).$$

Joshi and S\`a Barreto showed that for $w, w' \notin \partial X,$ the restriction of $U(t, w, w')$ to the diagonal is well defined.  Letting
$$i_{\Delta} : \R \times \Delta \to \R \times X \times X,$$
$$(t,w) \mapsto (t, w, w),$$
the pull back $i^* _{\Delta}$ is a Fourier integral operator of order $\frac{n}{4}.$  H\"ormander's transversal composition theorem \cite{ho} shows that $i^* _{\Delta} U(t)$ is a well defined distribution in $\R \times (X  - \partial X).$  Moreover, 
$$WF(i^* _{\Delta} U(t)) \subset \{ ((t, \tau), (w, \zeta - \eta)): \tau + q(w, \zeta) = 0, (w, \zeta) = \chi_t (w, \eta) \}.$$
For $\epsilon >0,$ let $X_{\epsilon} := \{x > \epsilon \}$ where $x$ is a boundary defining function for $\partial X.$  Let 
$$\pi : \R \times X_{\epsilon} \to \R,$$
$$(t, w) \mapsto t.$$
Integration over $w$ is equal to the push forward $\pi_*,$ so it is a Fourier integral operator.   By H\"ormander's theorem, 
$$T_{\epsilon} (t) := \int_{x> \epsilon} U(t, w, w) = \pi_* (i^* _{\Delta} U(t))$$
is a well defined distribution and in particular, \cite{JSB2} proved the following.

\begin{thm41} \em For $\epsilon > 0,$ the singular support of $T_{\epsilon}$ is contained in the set of periods of closed geodesics in $X_{\epsilon}.$  Moreover, there exists $\epsilon_0 > 0$ such that all closed geodesics of $(X, g)$ with period greater than zero are contained in $X_{\epsilon_0}.$  \em 
\end{thm41}

A further consequence of the iterative construction of the wave group in \cite{JSB2} is the existence of constants $C_j, j=1, \ldots, n-1,$ such that the limit
\begin{equation} \label{0trcos} \textrm{0-tr} (U(t)) = \lim_{\epsilon \to 0} \left[ \int_{x > \epsilon} U(t, w, w) - \sum_{j=1} ^{n-2} C_j \epsilon^{-j} + C_0 \log{\epsilon} \right] \end{equation}
exists.  This clearly depends on the choice of boundary defining function $x,$ but it does give a natural regularization of the trace $U(t)$ which is related to the length spectrum by the following.  
\begin{thm42}  \em The singular support of 0-tr $(U(t))$ is contained in the set of periods of closed geodesics of $(X, g).$   \em 
\end{thm42} 

Joshi and S\`a Barreto's Theorems 4.1 and 4.2 are the reason for the brevity of this note; the theorems allow us to generalize local dynamical arguments and results for compact manifolds to the asymptotically hyperbolic setting.

\section{Dynamical trace formula}

We recall some definitions from dynamics.  Let $SX$ denote the unit tangent bundle, and let $G^t$ be the geodesic flow on $SX.$  Since $X$ is complete and has negative sectional curvatures, the geodesic flow is Anosov \cite{anosov}; see for example \cite{bolton}, \cite{E}, \cite{E1}, \cite{K}.  Consequently, for each $\xi \in SX,$ $T(SX)_{\xi}$ splits into a direct sum 
\begin{equation} \label{split} T(SX)_{\xi} = E^s _{\xi} \oplus E^u _{\xi} \oplus E_{\xi},\end{equation}
where $E^s _{\xi}$ is exponentially contracting, $E^u _{\xi}$ is exponentially expanding, and $E_{\xi}$ is the one dimensional subspace tangent to the flow.  The \em Sinai-Ruelle-Bowen potential \em is a H\"older continuous function defined by 
\begin{equation} \label{srb} H(\xi) : = \frac{d}{dt} \big|_{t=0} \ln \det dG^t |_{E^u _{\xi}}.\end{equation} 
This potential is the instantaneous rate of expansion at $\xi.$  The topological pressure $\fp$ of a function $f: SX \to \R$ is defined as follows.  For large $T$ and small $\delta > 0,$ a finite set $Y \subset SX$ is $(T, \delta)$ separated if, given $\xi, \xi' \in SX, \xi \neq \xi',$ there is $t \in [0, T]$ with $d(G^t \xi, G^t \xi') \geq \delta.$  Here the distance on $SX$ is given by the Sasaki metric.  Then, 
\begin{equation} \label{press} \mathfrak{p}(f) = \lim_{\delta \to 0} \limsup_{T \to \infty}  T^{-1} \log \sup \left\{ \sum_{\xi \in Y} \exp \int_0 ^T f(G^t \xi) dt; \, Y \textrm{ is $(T, \delta)$ separated} \right\}.\end{equation}  
In the compact setting, the topological pressure of a function $f: SX \to \R$ may be equivalently defined by,   
$$ \mathfrak{p}(f) = \sup_{\mu} \left( h_{\mu} + \int f d\mu \right), $$
where the supremum is taken over all $G^t$ invariant measures $\mu,$ and $h_{\mu}$ denotes the \em measure theoretical entropy \em of the geodesic flow, \cite{Bo}, \cite{kh}.

The pressure of a function is a concept in dynamical systems arising from statistical mechanics which measures the growth rate of the number of separated orbits weighted according to the values of $f$ \cite{walt}.  In particular, $\mathfrak{p}(0) = h,$ the \em topological entropy \em of the geodesic flow,
\begin{equation} \label{topent} h = \lim_{\delta \to 0} \limsup_{T \to \infty}  T^{-1} \log \sup \# \{ Y \subset SX: Y \textrm{ is $(T, \delta)$ separated} \}.\end{equation}  

Note that for convex cocompact hyperbolic manifolds $\bbH^{n+1} / \Gamma,$ $h=\delta$ the exponent of convergence for the Poincar\'e series for $\Gamma$ which is also equal to the dimension of the limit set of $\Gamma.$  For the Sinai-Ruelle-Bowen potential $\mathfrak{p}(-H) = 0$ and the corresponding equilibrium measure which attains the supremum is the Liouville measure $\mu_L$ on the unit tangent bundle so that 
$$ h_{\mu_L} = \int_{SX} H d \mu_L.$$

\subsection{Two preliminary lemmas}

The local arguments of \cite{DG} together with Theorem 4.2 of \cite{JSB2} provide the leading terms in the renormalized wave trace, however, to bound the remainder term we adapt the local techniques of \cite{JPT} and \cite{JPS}, and this requires the following two lemmas.  

\begin{lemma}  \label{le:lyapunov} Let $(X, g)$ be a smooth, complete, $n+1$ dimensional Riemannian manifold whose sectional curvatures $\kappa$ satisfy 
$$-k_1^2 \leq \kappa \leq - k_2 ^2$$
for some $0 < k_2 \leq k_1.$  Then, the Poincar\'e map about a closed orbit $\gamma$ of the geodesic flow has eigenvalues $\lambda_i,$ $i=1, \ldots, 2n$ for which
$$e^{k_2 l(\gamma)} \leq |\lambda_i| \leq  e^{k_1 l(\gamma)}  \textrm{ for } i=1, \ldots , n
$$ and 
$$e^{- k_1 l(\gamma)} \leq |\lambda_i | \leq e^{-k_2 l(\gamma)} \textrm{ for } i = n+1 , \ldots, 2n,$$
where $l(\gamma)$ is the period (or length) of $\gamma.$  
\end{lemma} 

\textbf{Proof:  } Let $P_{\gamma}$ be the Poincar\'e map about the closed orbit $\gamma$ of the geodesic flow.  Since the flow is Anosov \cite{anosov}, $P _{\gamma}$ has $n$ expanding eigenvalues $\{\lambda_i\}_{i=1}^n$ and $n$ contracting eigenvalues $\{\lambda_i\}_{i=n+1} ^{2n}.$  We proceed to estimate the expanding eigenvalues using Rauch's Comparison Theorem.  Let $M_i$ be complete manifolds of dimension $n+1$ with constant negative curvature $-k_i^2$ for $i=1,2.$  Consider Jacobi fields $J$ and $J_i$ along any geodesic $\gamma_0$ on $X$ and $\gamma_i$ on $M_i$ such that 
$$J(0) = J_i (0) = 0, \quad \langle J'(0), \gamma_0'(0) \rangle = \langle J_i '(0), \gamma_i '(0) \rangle,$$
$$|J'(0)| = |J_i '(0)|.$$
Assume that $\gamma_0, \gamma_i$ do not have conjugate points on $(0,a]$ for some $a>0.$  Then, by Rauch's Comparison Theorem \cite{docarmo} for $t \in (0,a],$ 
\begin{equation} \label{eq:rauch} |J_2(t)| \leq |J(t)| \leq |J_1(t)|. \end{equation} 
By definition of the Lyapunov exponents \cite{pesin} for the Poincar\'e map on the geodesic flow, for the constant curvature manifolds $M_i,$ the expanding eigenvalues $\lambda_i$ of the Poincar\'e map about a closed geodesic $\tilde{\gamma}$ satisfy 
\begin{equation} \label{eq:cclya} | \lambda_i | = e^{k_i l(\tilde{\gamma})}. \end{equation}
Since the eigenvalues of the Poincar\'e map are determined by the Jacobi fields along the closed geodesics, it follows from (\ref{eq:rauch}) and (\ref{eq:cclya}) that the expanding eigenvalues for $P_{\gamma}$ satisfy
$$e^{k_2 l(\gamma)} \leq | \lambda_i | \leq e^{k_1 l(\gamma)}, \quad i=1, \ldots, n.$$
For each contracting eigenvalue $\lambda_i$ with $i \in \{n+1, \ldots, 2n\}$ there is an expanding eigenvalue $\lambda_{j(i)}$ with $j(i) \in \{1, \ldots, n\}$ such that 
$$|\lambda_i|^{-1} = |\lambda_{j(i)}|.$$
The inequality for the contracting eigenvalues follows immediately and completes the proof of the lemma.  \qed

The next result allows us to estimate the remainder in the trace formula by separating the periodic orbits and applying the stationary phase method of \cite{JPT} and \cite{JPS}.  This separation lemma is a generalization of Lemma 2.2.1 of \cite{JPT} to our $n+1$ dimension asymptotically hyperbolic variable negative curvature setting.  
\begin{lemma}  \label{le:sep}
Let $(X, g)$ be an asymptotically hyperbolic $n+1$ dimensional manifold with negative sectional curvatures.  Let $\cN (\gamma, \epsilon)$ denote the $\epsilon$ neighborhood of a geodesic $\gamma$ in the unit tangent bundle $SX$ with respect to the Sasaki metric.  Then there exist positive constants $T_0, B, $ and $\delta$ (depending only on the injectivity radius inj(X) and the curvature bounds) such that for any $T>T_0$ the sets $N(\gamma, e^{-BT})$ are disjoint for all pairs of closed geodesics $\gamma$ on $X$ with length $l_{\gamma} \in [T-\delta, T].$  
\end{lemma} 

\textbf{Proof:}
Since the sectional curvatures of any asymptotically hyperbolic manifold approach $-1$ at $\partial X,$ there exist $0 < k_2 \leq 1 \leq k_1$ so that 
\begin{equation} \label{curv} -k_1 ^2 \leq \kappa \leq -k_2 ^2, \end{equation}
for all sectional curvatures $\kappa$ on $X.$  Let $B > 2 k_1$ and choose $0< \delta < \textrm{inj}(X)/3,$ and let $T_0$ be such that $2e^{-k_1 T_0} < \delta.$  We proceed by contradiction.  For a given geodesic $\gamma_1,$ assume there is a second geodesic $\gamma_2$ with $T- \delta < l(\gamma_1) \leq l(\gamma_2) \leq T$ such that the corresponding neighborhoods intersect.  Assume the geodesics are not inverses of each other; by the choice of $\delta$ they cannot be integer multiples of each other unless they are inverses.  Let $\gamma_j (t), 0 \leq t \leq l(\gamma_j)$ denote the geodesic on $X,$ and let its corresponding lift to $SX$ be denoted by $\tilde{\gamma_j} (t) = (\gamma_j (t), \gamma_j '(t)).$  We may assume without any loss of generality that
$$ d_{SX} (\tilde{\gamma_2} (0), \tilde{\gamma_1} (0)) \leq 2 e^{-2 k_1 T}.  $$
For any $0 \leq t \leq l(\gamma_2),$ by the curvature assumptions and Lemma \ref{le:lyapunov}, 
\begin{equation} \label{eq:sx} d_{SX} (\tilde{\gamma_2}(t), \tilde{\gamma_1}(t) ) = d_{SX} (G^t \tilde{\gamma_2}(0), G^t (\tilde{\gamma_1}(0)) \leq 2 e^{-2k_1 T} e^{k_1 t} \leq 2e^{-k_1 T}, \end{equation}
where $G^t$ is the geodesic flow.  This implies
\begin{equation} \label{eq:dsx} d_X (\gamma_2 (t), \gamma_1 (t) ) \leq 2e^{-2k_1 T}. \end{equation} 
Consequently, the entire geodesics $\gamma_i$ lie in the $2e^{-k_1 T}$ neighborhood of each other.  Now, reparametrize $\gamma_1$ defining 
$$\beta_1 (s) := \gamma_1 (l_1 s/l_2), \quad 0 \leq s \leq l_2,$$
where $\gamma_i: [0, l_i] \to X.$  By the triangle inequality,
$$d(\gamma_2 (t) , \beta_1 (t) ) \leq d (\gamma_2 (t) , \gamma_1 (t)) + d(\gamma_1 (t), \beta_1 (t)) \leq$$
$$2e^{-k_1 t} + t \left( 1 - \frac{l_1}{_2} \right) \leq 2e^{-k_1 T} + \delta < \frac{2 \textrm{ inj}(X)}{3}.$$
For any $0 \leq t \leq l_2$ there exists a unique shortest geodesic $\alpha_t (s)$ in $X$ connecting $\gamma_2 (t)$ and $\beta_1 (t).$  Let the parameter $s \in [0,1]$ so that $\alpha_t (0) = \gamma_2 (t)$ and $\alpha_t (1) = \beta_1 (t).$  Define the mapping $\Phi(t,s) : [0, l_2] \times [0,1] \to X$ by the formula 
$$\Phi(t,s) = \alpha_t (s).$$
We will derive a contradiction by showing that $\Phi$ defines a homotopy between $\gamma_2(t)$ and $\beta_1 (t).$  First, $\Phi (t,0) = \gamma_2 (t),$ and $\Phi (t,1) = \beta_1 (t).$  Moreover, since both $\gamma_2$ and $\beta_1$ have period $l_2,$ we have 
$$\alpha_0 (s) = \alpha_{l_2} (s), \quad \forall \, s \in [0,1],$$
so that $\Phi (\cdot, s) $ is a closed curve in $X.$  Finally, $\Phi(t,s)$ is continuous since the function $d(\gamma_2 (t), \beta_1 (t))$ is a continuous function of $t.$  This shows that $\Phi$ is indeed a homotopy between $\gamma_2 (t)$ and $\beta_1 (t).$  Since $\beta_1$ is just a reparametrization of $\gamma_1,$ this shows that $\gamma_i$ lie in the same free homotopy class, which contradicts the fact that on a manifold with pinched negative curvature there is a unique closed geodesic in each free homotopy class  \cite{ehs}.  This contradiction proves the lemma.  
\qed

With these preliminary lemmas, we may now prove the dynamical wave trace formula.  

\subsection{Proof of Theorem \ref{theorem1}:}
Note that dependence of the renormalized trace of the wave group on the choice of boundary defining function is absorbed by the remainder term in the right side of the formula.  Fix $\epsilon > 0,$ and let $\phi \in \cC^{\infty} _0 ([\epsilon, \infty)).$  By \cite{JSB2} Theorem 4.2, the singular support of 0-tr $U(t)$ lies in $\{ k l(\gamma) : \gamma \in \cL_p\}.$  Since the arguments of \cite{DG} Theorem 4.5 are local, 0-tr $U(t)$ has an expansion at each singularity $T \in \{ k l(\gamma) : \gamma \in \cL_p\}$ with leading term 
\begin{equation} \label{dglead} \frac{ l(\gamma) \delta(|t| - k l(\gamma))}{ \sqrt{ | \det(I - \cP^k _\gamma ) | }}, \end{equation}
and moreover, 
$$0 \textrm{-}tr U(t) - \left( \sum_{\gamma \in \cL_p, \, k \in \N} \frac{ l(\gamma) \delta(|t| - k l(\gamma))}{ \sqrt{ | \det(I - \cP^k _\gamma ) |}} \right) = A(t)$$
is smooth for $t > 0.$

By the assumption of negative sectional curvatures together with the calculation of \cite{MMAH} which shows that all sectional curvatures approach $-1$ at $\partial X,$ there exist $0 < k_2 \leq 1 \leq k_1 $ so that
$$-k_1 ^2 \leq \kappa \leq -k_2^2,$$
for all sectional curvatures $\kappa.$  Therefore, the ``clean intersection'' condition of \cite{DG} is satisfied.  Since the arguments of \cite{DG} are localized to small neighborhoods around each $\gamma \in \cL,$ by (\ref{dglead}), \cite{DG}, and \cite{JSB2} Theorem 4.1, 
$$\int \phi(t) \textrm{0-tr} (U(t)) dt = \sum_{\gamma \in \cL_p, \, k \in \N, k l(\gamma) \in spt(\phi)} \frac{ l(\gamma) \phi(k l(\gamma))}{ \sqrt{ | \det(I - \cP^k _\gamma) | }} + A(\phi),$$
where $spt(\phi)$ is the support of $\phi$ and $A(\phi) = \int A(t) \phi(t) dt.$  By Lemma \ref{le:sep}, the periodic orbits are separated, and by \cite{JSB2} Theorem 4.1 the closed geodesics lie in a compact subset of $X,$ so we may apply the local estimates from the proof of Theorem 1.4.3 in \cite{JPT}, generalized to higher dimensions in \cite{JPS}.  Since $\phi$ has compact support in $[\epsilon, \infty)$ with $\epsilon > 0,$ by \cite{JPT} Proposition 3.2.1,
$$A(\phi) \leq C_{\phi},$$
where $C_{\phi}$ depends only on $||\phi||_{\infty}$ and $\epsilon.$  This completes the proof.  \qed

As a corollary to this theorem, we combine the dynamical trace formula with Borthwick's Poisson formula \cite{Borth} to produce a dynamical resonance-wave trace formula.  To state this result we recall a few definitions.  The Poisson formula relates the renormalized wave trace to the the poles, called \em resonances, \em of the meromorphically continued resolvent.  Closely related to the resolvent is the scattering operator whose poles essentially coincide with those of the resolvent; it is more convenient to state the trace formula in terms of scattering resonances.  Recall that the multiplicities of the resonances are given by 
$$m (\zeta) =   \textrm{ rank Res}_{\zeta} (\Delta - s(n-s))^{-1},$$
where $\zeta$ is a pole of the resolvent, $(\Delta - s(n-s))^{-1}.$  

\begin{defn} \label{scat} \em 
Let $(X, g)$ be an asymptotically hyperbolic $n+1$ dimensional manifold with boundary defining function $x.$  For $\mathfrak{R} s = \frac{n}{2}, \, s \neq \frac{n}{2},$ a function $f_1 \in \cC^{\infty} (\partial X)$ determines a unique solution $u$ of 
$$(\Delta - s(n-s) ) u = 0, \quad u \sim x^{n-s} f_1 + x^s f_2, \textrm{ as } x \to 0,$$
where $f_2 \in \cC^{\infty} ( \partial X).$  This defines the map called the \em scattering operator, \em  $S(s) : f_1 \mapsto f_2.$   \em 
\end{defn}

Heuristically, the scattering operator which is classically a scattering matrix, acts as a Dirichlet to Neumann map, and physically it describes the scattering behavior of particles.  The scattering operator extends meromorphically to $s \in \C$ as a family of pseudodifferential operators of order $2s - n.$  Renormalizing the scattering operator as follows gives a meromorphic family of Fredholm operators with poles of finite rank.  
$$\tilde{S} (s) : = \frac{\Gamma( s - \frac{n}{2})}{\Gamma(\frac{n}{2} - s) } \Lambda^{n/2 - s} S(s) \Lambda^{n/2 - s},$$
where 
$$\Lambda : = \frac{1}{2} (\Delta_h + 1)^{1/2}.$$
Above, $\Delta_h$ is the Laplacian on $\partial X$ for the metric $h(x) \big|_{x=0}.$  Note that this definition depends on the boundary defining function.  The multiplicity of a pole or zero of $S(s)$ is defined to be
$$\nu (\zeta) = - \tr [ \textrm{Res}_{\zeta} \tilde{S}'(s) \tilde{S}(s)^{-1} ].$$
The scattering multiplicities are related to the resonance multiplicities by 
$$\nu(\zeta) = m(\zeta) - m(n-\zeta) + \sum_{k \in \N} \left( \chi_{n/2 - k} (s) - \chi_{n/2 + k} (s) \right) d_k,$$
where 
$$d_k = \dim \ker \tilde{S} (\frac{n}{2} + k),$$
and $\chi_{p}$ denotes the characteristic function of the set $\{p\}.$  Recall the boundary $(\partial X, [h_0])$ does not have a unique metric but rather admits a conformal class $[h_0]$ of metrics determined by $h(x) \big| _{x=0}$ for a totally geodesic boundary defining function $x.$  Graham, Jenne, Manson, and Sparling \cite{gjms} defined natural conformally invariant powers of the Laplacian on conformal manifolds.  Using scattering theory, Guillop\'e and Zworski \cite{gzscat} identified these operators to poles of the scattering operator on asymptotically hyperbolic manifolds.  Let $P_k$ denote the $k^{th}$ conformal power of the Laplacian on the conformal infinity $(\partial X, [h_0]).$  When the dimension is odd, these are not naturally defined; see \cite{GN} for their definition.  The numbers $d_k$ above are equivalently given by the dimensions of the kernels of these conformal operators $P_k$ acting on the conformal infinity $(\partial X, [h_0]).$  

The set of resolvent resonances will be denoted $\cR,$ while the set of scattering resonances 
$$\cR^{sc} : = \cR \bigcup_{n=1} ^{\infty} \{ \frac{n}{2} - k \textrm{ with multiplicity } d_k \}.$$  Finally, we recall that for a discrete torsion free group $\Gamma$ of isometries of $\bbH^{n+1},$ the quotient $\bbH^{n+1} / \Gamma$ is said to be \em convex cocompact \em when its convex core is compact; a compact perturbation is defined as follows.  

\begin{defn} \label{pert} \em An asymptotically hyperbolic manifold $(X, g)$ is a \em compact perturbation \em of a convex cocompact hyperbolic manifold if there exists a convex cocompact manifold $(X_0, g_0)$ (possibly disconnected) such that $(X-K, g) \cong (X_0 - K_0, g_0)$ for some compact sets $K \subset X$ and $K_0 \subset X_0.$  \em 
\end{defn}

\subsection{Proof of Corollary \ref{corollary1}:}
This corollary is an immediate consequence of Theorem 1.2 of \cite{Borth} and our dynamical trace formula.  We note that in terms of the resolvent resonances the spectral side of the trace formula becomes 
\begin{equation} \label{eq:res} \frac{1}{2} \sum_{s \in \cR } m(s) e^{(s - n/2)|t|} + \frac{1}{2} \sum_{k \in \N} d_k e^{-k|t|} + B(t).\end{equation} \qed  

One application of this trace formula is a lower bound for the length spectrum counting function in the following corollary.  This lower bound is not optimal; we show in the prime orbit theorem that the length spectrum counting function grows like $e^{hT}/hT,$ where $h$ is the topological entropy of the geodesic flow (\ref{topent}).  However, we include this corollary for its proof which relies on the exponential growth in the spectral side of the trace formula to estimate the length spectrum counting function and thereby distinguishes hyperbolic scattering from Euclidean scattering, since in Euclidean scattering the spectral side of the trace formula is purely oscillatory.  

\subsection{Proof of Corollary \ref{corollary2}:}  
Let $\phi \in \cC^{\infty} _0 (\R^+)$ and $T >>0$ so that the support of $\phi$ lies in $[T-2, T+1]$ and $\phi \equiv 1$ on $[T-1, T].$  Using $\phi$ as a test function in the trace formula (\ref{eq:res}) the spectral side gives 
$$\frac{1}{2} \sum_{ s \in \cR} \int e^{(s- n/2)t} \phi(t) d t + \frac{1}{2} \sum_{k \geq 1} \int e^{-d_k t} \phi(t) d t + O(1).$$  
The dominant term is asymptotic to 
\begin{equation} \label{eq:dom} \frac{1}{2} e^{(s_0 - \frac{n}{2} ) T}, \quad T >>0, \end{equation} 
and the remainder is 
\begin{equation} \label{eq:rem} O(e^{\epsilon T}) \textrm{ for some fixed $\epsilon \in (0, s_0 - n/2)$ which depends on $\sigma_{pp} (\Delta).$} \end{equation}
In particular, if 
$$\sigma_{pp} (\Delta) = 0 < \Lambda_0 < \Lambda_1 < \ldots \frac{n^2}{4},$$
then the remainder is $O(e^{(s_1 - \frac{n}{2})T})$ where $\Lambda_1 = s_1 (n-s_1),$ and $n/2 < s_1 < s_0.$   If 
$$\sigma_{pp} (\Delta) = \{\Lambda_0\},$$
then by Theorem 1.1 of \cite{Borth}, the remainder is $O(e^{\epsilon T})$ for any $ \epsilon > 0.$  The dynamical side of the trace formula with test function $\phi$ is 
\begin{equation} \label{eq:cortrace} \frac{1}{2} \sum_{\gamma \in \cL_p,\,  k \in \N \, :  \, k l(\gamma) \in (T-2, T+1) } \frac{l(\gamma) \phi(k l(\gamma) ) }{\sqrt{|\det( I - \cP^k_{\gamma} )|}} + \cO(1). \end{equation} 
Since $X$ is a negatively curved compact perturbation of a hyperbolic manifold, we may apply Proposition \ref{le:lyapunov} to estimate $| \det(I - P^k _{\gamma} )|.$ 
Since $\cP^k _{\gamma}$ has expanding eigenvalues $\lambda_1, \ldots, \lambda_n$ and contracting eigenvalues $\lambda_{n+1}, \ldots, \lambda_{2n},$ 
$$|\det(I - \cP^k _{\gamma})| = \prod_{1} ^{2n} |1 - \lambda_i| = \prod_{i=1} ^n |\lambda_i| \prod_{j=1} ^n \left|1 - \frac{1}{|\lambda_j|}\right| \prod_{k=n+1} ^{2n} |1 - \lambda_k|.$$
By Lemma \ref{le:lyapunov}, 
\begin{equation} \label{eq:pupper} |\det(I - \cP^k _{\gamma})| \leq e^{n k_1 k l(\gamma)} \left(1- O(e^{-k_2 k l(\gamma)} ) \right) \end{equation} 
and
\begin{equation} \label{eq:plower}|\det(I- \cP^k _{\gamma} ) | \geq e^{n k_2 k l(\gamma)} \left( 1 - O(e^{-k_2 k l(\gamma)}) \right).\end{equation} 
Let 
$$N(T) = \# \{ \gamma \in \cL : l(\gamma) \leq T \}.$$
Then, (\ref{eq:dom}), (\ref{eq:rem}), (\ref{eq:cortrace}), (\ref{eq:pupper}), and (\ref{eq:plower}) imply the growth estimates for $N(T),$
$$\limsup_{T \to \infty} \frac{T N(T)}{e^{(s_0 + (k_1 - 1)n/2)T}} \leq 1 \leq \liminf_{T \to \infty} \frac{T N(T)}{e^{(s_0 + (k_2 - 1)n/2)T}}.$$ \qed

\textbf{Remarks:  }
Note that for convex cocompact hyperbolic manifolds, $k_2 = k_1 = 1,$ and it is known that when the topological entropy $h > n/2,$ $s_0 = h$ \cite{N}, so we recover the prime orbit theorem of \cite{PLength}.  We note that when $\sigma_{pp} (\Delta) = \emptyset,$ the corollary is vacuous since the spectral side of the trace formula is no longer dominated by the term $e^{(s_0 - n/2)t}.$  However, once we have proven the prime orbit theorem, we will show that the existence of pure point spectrum of the Laplacian is determined by the topological entropy of the geodesic flow and the curvature buonds.  

Another application of the trace formula is to counting resonances in strips.  See for example Theorem 1.3 of \cite{P} and Theorem 2 of \cite{GZ} which we may now extend to compact perturbations of convex cocompact negatively curved hyperbolic manifolds in arbitrary dimension.  Those results use the existence of only \em one \em closed geodesic; by incorporating more of the length spectrum, one expects these estimates to be improved to give a fractal Weyl law with exponent determined by the entropy of the geodesic flow.  This remains an interesting open problem.  


\section{Dynamical Zeta Function}

The dynamical zeta function is to the geodesic length spectrum as the Riemann zeta function is to the prime numbers.  In particular, let
\begin{equation} \label{dz} Z(s) = \exp  \left( \sum_{\gamma \in L_p} \sum_{k \in \N} \frac{e^{-k s l_p(\gamma)}}{k} \right), \end{equation}
where $L_p$ consists of primitive closed orbits of the geodesic flow and $l_p(\gamma)$ is the primitive period (or length) of $\gamma \in L_p.$  This definition is the same as Parry and Pollicott's dynamical zeta function for Axiom A flows, \cite{POT}.  It is also interesting to consider the following weighted dynamical zeta function, 
$$ \tilde{Z} (s) = \exp \left( \sum_{\gamma \in L_p} \sum_{k \in \N} \frac{e^{-k s l_p(\gamma)}}{k \sqrt{| \det(I - \cP^k _{\gamma}) |}} \right),$$
where $\cP^k _{\gamma}$ is the $k$-times Poincar\'e map of the geodesic flow around the closed orbit $\gamma.$  The weighted zeta function is particularly interesting for its connections to the resonances of the resolvent; \cite{P} and \cite{GN} used the Hadamard factorization of this zeta function to prove Selberg trace formulae for convex cocompact hyperbolic manifolds.  

\subsection{Proof of Theorem \ref{theorem2}:} 
We adapt the methods of \cite{cm} to our setting.  Note that the sum is 
$$\sum_{\gamma \in L} k(\gamma)^{-1} \exp{(- \int_{\gamma} s dt)} = \sum_{r \in \N} a_r,$$
where $L$ consists of all closed orbits of the geodesic flow, $k(\gamma)$ is the multiplicity of $\gamma,$ and 
$$a_r = \sum_{\epsilon(r - 1/2) \leq l(\gamma) < \epsilon(r + 1/2)} k(\gamma)^{-1} \exp{(-\int_{\gamma} s dt)}.$$
The arguments of \cite{f} for Lemma 2.8 are entirely local and show that 
$$a_r \leq \frac{\exp{-(r \epsilon \mathfrak{p}(-s))}}{r \epsilon}.$$
Moreover, by \cite{walt}, $\mathfrak{p}(-s) = \fp(0) - s,$ so that the series converges absolutely when 
$$\mathfrak{R}(s) > \fp(0).$$
For the weighted dynamical zeta function, note that $\cP_{\gamma}$ has expanding eigenvalues $\lambda_1, \ldots, \lambda_n,$ and contracting eigenvalues $\lambda_{n+1}, \ldots, \lambda_{2n},$ and
$$|\det(I - \cP_\gamma)| = \prod_{1} ^{2n} |1 - \lambda_i| = \prod_{i=1} ^n |\lambda_i| \prod_{j=1} ^n \left|1 - \frac{1}{|\lambda_j|} \right| \prod_{k=n+1} ^{2n} |1 - \lambda_k|.$$
By Lemma \ref{le:lyapunov}, 
$$|\lambda_i|^{-1} \leq e^{- k_2 l(\gamma)}, \quad i=1, \ldots, n,$$
and 
$$|\lambda_i| \leq e^{-k_2 l(\gamma)}, \quad i=n+1, \ldots 2n.$$
Therefore, 
$$ \lim_{|l(\gamma)| \to \infty} \frac{\prod_{i=1} ^n |\lambda_i|}{ |\det(I - \cP_\gamma)| } = 1,$$
so we may replace $|\det(I - \cP_\gamma)|$ by this product of expanding eigenvalues.  Since $H$ is the rate of expansion of volume in $E^u,$ the summand for $\gamma$ is
$$k(\gamma)^{-1} \exp( \int_{\gamma} \frac{1}{2} H - s ).$$
Then, similarly defining $a_r,$ Lemma 2.8 of \cite{f} shows that 
$$a_r \leq \frac{\exp( r \epsilon \fp(\frac{1}{2} H - s) )}{r \epsilon},$$
so that the series converges absolutely when $\fp(\frac{1}{2} H - s) < 0.$  Since $\fp(\frac{1}{2} H - s) = \fp(\frac{1}{2} H) - s,$ this shows that the series converges absolutely when 
$$\mathfrak{R}(s) > \fp(\frac{1}{2} H).$$ 
\qed

\textbf{Remark:}  As observed in \cite{cm}, if the exponent $\frac{1}{2}$ in the denominator of the weighted dynamical zeta function  is replaced by $t \in \R,$ then by the preceding arguments $\tilde{Z}(s)$ converges absolutely for $\mathfrak{R}(s) > \fp(- t H).$  

\section{Prime Orbit Theorem}

To prove the prime orbit theorem we require further definitions to describe the geodesic flow.  The following definitions are from \cite{E1} and \cite{E2}; see also \cite{BO} and \cite{ehs}.  For $\xi \in SX,$ recall the \em positive prolongational limit set, \em 
$$P^+ (\xi) = \{ y \in SX: \textrm{ for any neighborhoods $O, U$ of $\xi, y,$ respectively, there is a}$$
$$\textrm{sequence } t_n \subset \R, t_n \to \infty, \textrm{ such that } G^{t_n} (O) \cap U \neq \empty \}.$$
Then $\xi$ is non-wandering if $\xi \in P^+ (\xi).$  The flow is \em topologically transitive \em on $\Omega \subset SM$ if for any open sets $O, U \subset \Omega,$ there exists $t \in \R$ such that $G^t(U) \cap O \neq \emptyset.$  The flow is \em topologically mixing \em if there exists $A > 0$ so that for all $|t| > A,$ $G^t (U) \cap O \neq \emptyset.$  

A closed invariant set $\Omega \subset SX$ without fixed points is \em hyperbolic \em if the tangent bundle restricted to $\Omega$ is a Whitney sum 
$$T_{\Omega} SX = E + E^s + E^u$$
of three $TG^t$ invariant sub-bundles, where $E$ is the one dimensional bundle tangent to the flow, and $E^s,$ $E^u$ are exponentially contracting and expanding, respectively:  
$$||TG^t (v) || \leq K e^{- \lambda t} ||v|| \textrm{ for } v \in E^s, t \geq 0,$$
$$||T G^{-t} (v) || \leq K e^{- \lambda t} ||v|| \textrm{ for } v \in E^u, t \geq 0.$$
In \cite{POT}, Parry and Pollicott defined a \em basic set \em to be a topologically transitive hyperbolic set $U$ with no fixed points for which periodic orbits are dense and which admits an open set $O \supset U$ so that $U = \cap_{t \in \R} G^t O.$  

\subsection{Proof of Theorem \ref{theorem3}:}  With the work of Bishop-O'Neill \cite{BO}, Eberlein \cite{E}, \cite{E1}, \cite{E2}, and Eberlein-O'Neill \cite{eo} on ``visibility manifolds'' (complete manifolds with non-positive curvature), we are able to give a quick proof of the prime orbit theorem.  For $X$ asymptotically hyperbolic and negatively curved, the non-wandering set $\Omega \subset SX$ is closed and invariant under the flow; see \cite{E1} page 502.  Moreover, by Theorems 3.9 and 3.10 in \cite{E1}, the periodic vectors are dense in $\Omega.$\footnote{In the cases where either $X$ has no closed geodesics or only one primitive closed geodesic, the theorem is vacuous.}  By \cite{E1} Theorem 3.13, $\Omega$ is connected, and by \cite{E1} Theorem  3.11, the geodesic flow restricted to $\Omega$ is topologically transitive.  Since the flow is Anosov, $\Omega$ is a hyperbolic set (\ref{split}).  Since the periodic orbits are dense in $\Omega,$ and $\Omega$ is closed, by Theorem 4.1 of \cite{JSB2}, $\Omega$ is a compact subset of $SX.$  Cleary $\Omega$ may not have fixed points for the geodesic flow, and by definition of the non-wandering set, for any open neighborhood $O \supset \Omega,$ $\cap_{t \in \R} G^t O = \Omega.$  Consequently, $\Omega$ is a basic set.  By \cite{POT} Proposition 1, the flow restricted to $\Omega$ is topologically mixing.  Since all Anosov flows are \`a priori Axiom A flows, the geodesic flow restricted to $\Omega$ is a topologically mixing Axiom A flow restricted to a basic set and satisfies the hypotheses of \cite{POT} Theorems 1 and 2.  Applying the results of Theorems 1 and 2 from \cite{POT} completes the proof of our theorem.  \qed

\subsection{Proof of Corollary \ref{corollary3}:}  
The proof consists of applying the prime orbit theorem and analyzing the dominant terms on the spectral and dynamical sides of the trace formula in Corollary \ref{corollary1}.  First, assume the topological entropy satisfies $h > \frac{n k_1}{2}.$  Considering a test function $\phi$ as in the proof of Corollary \ref{corollary2}, the dynamical side of the trace formula is bounded below by a constant multiple of 
$$e^{hT - nk_1 T/2} \quad \textrm{as } T \to \infty.$$
This implies that the spectral side of the trace formula must also have exponential growth.  Since the only terms with positive exponents come from $\sigma_{pp}(\Delta),$ there must exist $\Lambda_0 =s_0 (n-s_0) \in \sigma_{pp}(\Delta)$ with 
$$s_0 - n/2 \geq h - nk_1/2 \implies s_0 \geq h + \frac{n(1-k_1)}{2}.$$
This proves the first statement of the corollary.  Next, assume the topological entropy satisfies $h \leq \frac{n k_2}{2}.$  Considering the same test function $\phi,$ the dynamical side of the trace formula is bounded above by 
$$T e^{(h - nk_2/2)T},$$
so does not have exponential growth.  Therefore, the spectral side of the trace formula cannot have exponential growth.  Since any $s$ with $\Lambda = s(n-s) \in \sigma_{pp}(\Delta)$ is strictly greater than $n/2,$ we necessarily have $\sigma_{pp} (\Delta) = \emptyset$ which completes the proof.  \qed

\section{Concluding Remarks and Further Directions}  
This work was motived by \cite{GZ}, \cite{GN}, and \cite{G1}; our aim was to prove trace formulae for asymptotically hyperbolic $n+1$ dimensional manifolds and to understand the relationship between the resonances and dynamics of these manifolds.  Although our trace formulae do not allow an explicit expression for the remainder terms, nonetheless, like classical trace formulae they provide a deep and beautiful connection between the Laplace and length spectra with applications to computing remainder terms for both the length and resonance counting functions; see, for example \cite{JPT} and \cite{GN}.  Another application is to counting resonances in regions of $\C$ corresponding to physical phenomena; see \cite{GZ}.  It would be interesting to study the remainder terms in our formulae in greater depth; it is almost certain that an explicit formula for the remainder does not exist in this context of variable curvature, but perhaps one may show exponential decay in the remainder at infinity.  It would also be interesting to study the behavior as $t \to 0.$  Ideally, we would like to generalize our trace formulae to all asymptotically hyperbolic manifolds.  The dynamical trace requires hypotheses on the geodesic flow to allow summation of periodic orbits and to control the remainder term; assuming globally negative (but not necessarily constant) curvature guarantees this, but based on \cite{anosov}, \cite{E}, and \cite{K}, we expect a weaker hypothesis to suffice.  Such a hypothesis may be quite technical.  The resonance formula is more delicate and remains an open problem for asymptotically hyperbolic manifolds.  Borthwick's recent Poisson formula \cite{Borth} is progress in this direction, although Guillarmou's careful study of the scattering phase \cite{G1} indicates that the Poisson Formula for asymptotically hyperbolic manifolds is a subtle and elusive task.  We hope that our work is a useful contribution to understanding the Laplace and length spectra on asymptotically hyperbolic manifolds, and we encourage readers interested in pursuing the many open problems which remain for these spaces.

\end{document}